\documentclass[11pt]{gen-j-l}
\usepackage{amsmath}
\usepackage{amsfonts,amssymb,amsthm}
\usepackage{graphicx}
\newcommand{\half}{\frac12}
\newcommand{\crat}{\mathfrak{c}}
\newcommand{\thm}[2]{\begin{#1} #2 \end{#1}}

\def\card{\mathop {\rm card}\nolimits}
\newcommand{\QED}{\end{proof}}
\newcommand{\bd}{\textbf}
\newcommand{\sh}{\mathfrak{s}}

\newcommand{\hthree}{\mathbb{H}^3}
\newcommand{\reals}{\mathbb{R}}
\newcommand{\comps}{\mathbb{C}}
\newcommand{\ccirc}[1]{\mathcal{C}_{#1}}
\newcommand{\htwo}{\mathbb{H}^2}
\newcommand{\cbar}{\overline{\mathbb{C}}}

\begin{document}
\title{f Intrinsic Geometry of Convex Ideal Polyhedra in Hyperbolic
3-space} 
\copyrightinfo{1994}{Igor Rivin}
\date{}
\author{Igor Rivin}
\address{Mathematics department, University of Manchester, Oxford
Road, Manchester}

\thanks{}

\begin{abstract}
I describe a simple relationship between triangulations in $\mathbb{c}$ and
ideal polyhedra in $\mathbb{H}^3.$ 
I produce a complete intrinsic characterization of convex polyhedra in
hyperbolic 3-space $\hthree$ with all vertices on the sphere at
infinity. I also show that such polyhedra are uniquely determined by
their intrinsic metric. 
\end{abstract}

\maketitle

\section{Introduction}
\label{intro}

In this paper I study the intrinsic geometry of convex polyhedra in
three-dimensional hyperbolic space $\hthree$, with all vertices on the
sphere at infinity $S_\infty^2$. Such a polyhedron $P$ is homeomorphic to
the sphere $\mathbb{S}^2$ with a number of punctures (corresponding to the
vertices of $P$). It is not hard to see that $P$ is a complete 
hyperbolic surface of finite area.

In this paper I prove the following converse:

\thm{theorem}
{
\bd{Characterization of ideal polyhedra} \label{charRI}
Let $M_N$ be a complete hyperbolic surface of finite area,
homeomorphic to the $N$ times punctured sphere. Then $M_N$ can be
isometrically embedded in $H^3$ as a convex polyhedron $P_N$ with all
vertices on the sphere at infinity.
}

Furthermore,

\thm{theorem}
{
\bd{Uniqueness of realization} \label{Uniq}
The polyhedron $P_N$ promised by Theorem \ref{charRI} is unique, up to
congruence.
}

In fact, Theorem \ref{Uniq} is a special case of the following more
general result:

\thm{theorem}
{
\bd{Uniqueness of generalized polyhedra}\label{Uniq1}
A {\em generalized polyhedron} $P_N$ in $H^3$ is determined by its
intrinsic metric, up to congruence.
}

A \textit{generalized polyhedron} is a polyhedron some of whose vertices
are inside $\hthree$, some are on the ideal boundary of $\hthree$ and
some are beyond the ideal boundary. This is best visualized in the
projective model of $\hthree$ as a (Euclidean) polyhedron, all of
whose edges  intersect the unit ball $\hthree$. 
The derivation of Theorem \ref{Uniq1} from Theorem \ref{links} is
essentially the same as that of Theorem \ref{Uniq}.

The proof of Theorem \ref{charRI} uses the Invariance of Domain
Principle of A.~D.~Aleksandrov. This is as follows:

Given a map $f:A\rightarrow B$ between topological spaces $A$ and $B$,
then $f$ is onto, provided the following criteria are satisfied:

\begin{enumerate}

\item{} The image of $f$ is non-empty.

\item{} $f$ is continuous.

\item{} $f$ maps open sets in $A$ to open sets in $B$.

\item{} $f$ maps closed sets in $A$ to closed sets in $B$.

\item{} $B$ is connected.

\end{enumerate}

\begin{proof}[of Theorem \ref{charRI} (Outline)]
Let $M_\mu$ be a surface
homeomorphic to the $N$-times punctured 
2-sphere, with a {\em marking} $\mu$, that is, a labelling of
the punctures.
Let $\mathcal{P}^N$ be the space of convex ideal polyhedra
in $\hthree$, parametrized by the positions of their vertices on the
sphere at infinity of $\hthree$ (interpreted as the Riemann sphere
$\cbar$). Three of the vertices are fixed at $0$, $1$ and 
$\infty$, which eliminates the action of the isometry group of $\hthree$.
$\mathcal{P}^N$ is easily seen to be a $2 N - 6$ dimensional manifold.
$\mathcal{P}^N$ plays the role of $A$ in the invariance of domain
principle. We will abuse notation and view a polyhedron $P$ both as a
geometric object and as a polyhedral isometric embedding of $M$ into
$\hthree$. $P$ will inherit the labelling of vertices from $\mu$.

Let $\mathcal{T}^N$ be the set of complete, finite volume hyperbolic
structures on $M_\mu$.

This set is parametrized by
{\em shears} along the edges of a geodesic triangulation. This
parametrization is explained in detail in section \ref{sec-intrinsic}.
It will also be shown (Theorem \ref{teich}) that $\mathcal{T}^N$ is a $2
N - 6$ dimensional contractible manifold. Although many
of the results of section \ref{sec-intrinsic} are known to
Teichm\"uller theorists, they are so elementary in this
particular setting that it was impossible to resist including a full
exposition. 

 $\mathcal{T}^N$ will play the
role of $B$ in the invariance of domain principle.

The role of the map $f$ will be played by the map $\mathfrak{g}: \mathcal{P}^N
\rightarrow \mathcal{T}^N$. $\mathfrak{g}(P)$ is $P$ viewed as an abstract
Riemannian manifold. The continuity of $\mathfrak{g}$ with respect to the chosen
coordinate systems on $\mathcal{P}^N$ and $\mathcal{T}^N$ is the content of
Theorem \ref{mapcont}.

Since $\mathcal{P}^N$ and $\mathcal{T}^N$ are manifolds of the same
dimension and $\mathfrak{g}$ is continuous, Theorem \ref{Uniq}
shows that $\mathfrak{g}$ is an open map. That $\mathfrak{g}$ is closed is
the content of Theorem \ref{closed}. \end{proof}

The theory developed in Section \ref{sec-cross} is of independent
interest. In particular, it leads to a very simple derivation of a set
of conditions satisfied by dihedral angles of an ideal polyhedron
(Theorem 3.12).

\section{Hyperbolic geometry of the $N$-times punctured sphere}
\label{sec-intrinsic}
\subsection{Geometry of triangles}

Let $ABC$ be an ideal triangle in $\htwo$. Pick a point $p$ on the
geodesic $AB$. How far is $p$ from $A$? This question turns out to
make sense in $\htwo$:

Consider the unique horocycle $h_A$ centered at $A$ and passing
through $p$. $h_A$ will intersect $AC$ in a point $q$. Define
$\mathfrak{D}_{ABC}(p)$ to be the distance along $h_A$ between $p$ and $q$.
$\mathfrak{D}_{ABC}$ has the following important property:

\thm{Lemma}{\label{horos}
Let $p_1$ and $p_2$ be two points on $AB$. The hyperbolic
distance between $p_1$ and $p_2$ is equal to 
$|\log (\mathfrak{D}_{ABC}(p_1)/\mathfrak{D}_{ABC}(p_2))|.$
}

\begin{proof}
Use the upper half-space model of $\htwo$. Recall that the
hyperbolic metric is related to the Euclidean metric on the upper
half-space by $ds_h = |dz|/\Im z$. This means, in particular, that
the hyperbolic distance between $z_1 = x + i y_1$ and $z_2 = x + i
y_2$ is $|\log(y_1/y_2)|$.  

By a hyperbolic isometry $A$ can be sent to $\infty$, $B$ to $0$ and
$C$ to $1$. Horocycles centered on $A$ are then simply horizontal
lines, and if $p = x + i y$, then $\mathfrak{D}(p) = 1/y$, since the metric on
the horocycle around infinity through $p$ is then simply the standard
metric on the real line, rescaled by $1/y$. The assertion of the Lemma
now follows. \end{proof}

Lemma \ref{horos} gives a way to quantify the ways in which two
ideal triangles can be joined together along a side to form an ideal
quadrilateral. Intuitively, ideal triangles $ABC$ and $ADC$ can slide
with respect to each other along the common side $AC$. Pick a point
$p$ on $AC$. If $\mathfrak{D}_{ACB}(p) = \mathfrak{D}_{ACD}(p)$ ($\mathfrak{D}$ is
taken with respect to the vertex $A$), then we say that
$ABC$ and $ADC$ are joined without a shear (and it is easy to see
that reflection in $AC$ will send $B$ to $D$ and vice versa).
Otherwise, $ABC$ and $ADC$ are joined with shear $\log
(\mathfrak{D}_{ACB}(p)/\mathfrak{D}_{ACD}(p))$. It is clear that shear doesn't
depend on which of the vertices $A$ or $C$ is taken as the center of
the horocycles. Henceforth the shear between triangles $t_1$ and $t_2$ 
will be denoted by $\sh(t_1, t_2).$ By abuse of notation
$\sh(ABCD)=\sh(ABC, ADC)$.

All of the above is quite easily seen in the upper half-space model
of $\htwo$ -- see Figure 1. Figure 1 also demonstrates Lemma
\ref{crossr1}.

\thm{Lemma}{\label{crossr1}
If the shear is $\alpha$, and the triangle $ABC$ is
positioned so that $A = \infty$, $B=1$, $C=0$, then $\log |D|=
\alpha$.
}

\begin{figure}
\label{shearf}
\includegraphics*[width=3in,height=3in]{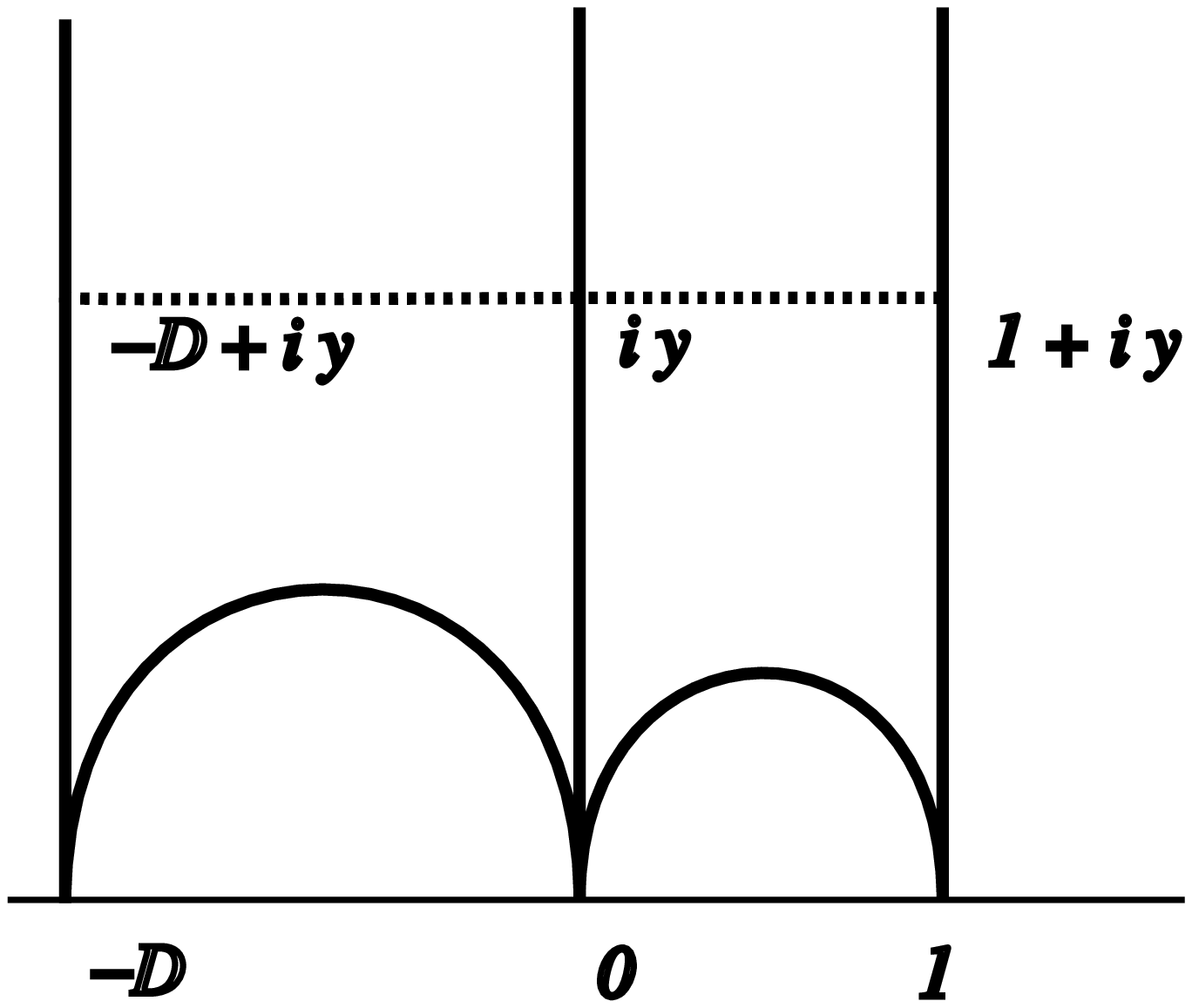}
\caption{Shear in the upper half-space model of $\htwo$}
\end{figure}

This also shows that the following fundamental lemma:

\thm{Lemma}{\label{crossr}
$\sh(ABCD)$ is equal to the log of the absolute value of the {\bf
cross ratio}  $[C, B, D, A]$.
}

Recall that cross-ratio of $z_1, z_2, z_3, z_4$ is defined to be
$$[z_1, z_2, z_3, z_4] = \frac{(z_1-z_3)(z_2 - z_4)}{(z_1-z_2)(z_3
- z_4)}$$

\begin{proof}
Both $\sh(ABCD)$ and $[C, B, D, A]$ are invariant under
the M\"obius group, and since any three points can be transformed to
$0$, $1$, and $\infty$ by a M\"obius transformation, the Lemma follows
from the obvious special case where $A=\infty$, $B=1$, and $C=0$. \end{proof}

\subsection{Geometry of triangulations}

Recall that $M_\mu$ is a surface homeomorphic to the 2-sphere with
$N$ punctures, together with its marking.

Let $T$ be a triangulation of $M_\mu$. Then the above discussion
serves to parametrize all the complete hyperbolic structures on $(M,
\mu)$ where the faces of $T$ are ideal triangles --- to each edge of
$T$ we associate the shear of the two abutting faces of $T$. This
information specifies the geometry completely. On the other hand, it
is not quite true that any assignment of real numbers to edges of $T$
corresponds to a complete hyperbolic structure on $M_\mu$ with those
numbers as shears -- it is necessary and sufficient that the sum of
shears around any cusp add up to zero (pick a horocycle $h$ centered
at a vertex $v$ of $T$. In order for a hyperbolic structure to be
complete, $h$ must close up). Thus if there are $V$ vertices and $E$
edges of $T$, the set of hyperbolic structures on $M_\mu$ such that
$T$ is an ideal triangulation is naturally parametrized as $R^{E-V}$.
Note that an Euler's formula computation yields  $E-V = 2 V - 6$,
so the dimension of the space of hyperbolic structures depends only on
the number of cusps.

The following lemma shows that the reliance on the particular
triangulation $T$ in the above discussion is not critical -- any other
topological ideal triangulation will do as well:

\thm{Lemma}{\bd{straightening} \label{strght}
Any topological triangulation $\tilde{T}$ with all vertices at cusps
of $M_\mu$
can be straightened to a geodesic triangulation. 
}

\begin{proof}
All that is necessary to show is that if $v_1$, $v_2$,
$v_3$ and $v_4$ are cusps and there are non-intersecting curves
$\gamma_1$ connecting $v_1$ and $v_2$ and $\gamma_2$ connecting $v_3$
and $v_4$, then the corresponding geodesics also don't intersect. That
this last statement is true can be observed by examining the geometry
of universal cover of $S$. The postulated curves $\gamma_1$ and
$\gamma_2$ exist if and only if the lifts of $v_1$ and $v_2$ do
not separate the lifts of $v_3$ and $v_4$ on the circle at
infinity of $H^2$. In that case, however, it is seen that the
corresponding geodesic segments do not intersect either. \end{proof}

Thus, every topological triangulation of $M_\mu$ with vertices at the
cusps ({\em topological ideal triangulation}) corresponds to a
coordinate system on the space of hyperbolic metrics on $S$ with
cusps at the prescribed vertices. 

\thm{Definition}
{
The space of complete hyperbolic structures on $M_\mu$ shall 
be denoted by $\mathcal{T}^N$. A {\em shear coordinate system
corresponding to a triangulation $T$} on $\mathcal{T}^N$ is the map
$\mathcal{C}_T: \mathcal{T}^N\rightarrow \reals^{2 N - 6}$ associating to a
particular metric its shears along the straightened edges of $T$. 

}

\thm{theorem}{ \label{ctris}
For any two triangulations $T_1$ and $T_2$, the map $c_{T_1 T_2} =
\mathcal{C}_{T_2} \circ \mathcal{C}_{T_1}^{-1}$ is a continuous function from
$\reals^{2 N - 6}$ to itself. 
}

To prove Theorem \ref{ctris} it will first be necessary to understand
{\em triangulation graph} of the sphere with $N$ vertices.

\thm{Definition}{\label{whitem}
Let $T$ be a triangulation, and $ABC$ and $ADC$ be two of the
triangles of $T$ sharing the edge $AC$. Then the {\em Whitehead move}
$w_{ABCD}$ transforms $T$ into a triangulation $T'$, where the
triangles $ABC$ and $ADC$ are replaced by triangles $BAD$ and $BCD$
(in other words the diagonal of the quadrilateral $ABCD$ is
``flipped'').
}

\begin{figure}
\includegraphics*[width=3in,height=3in]{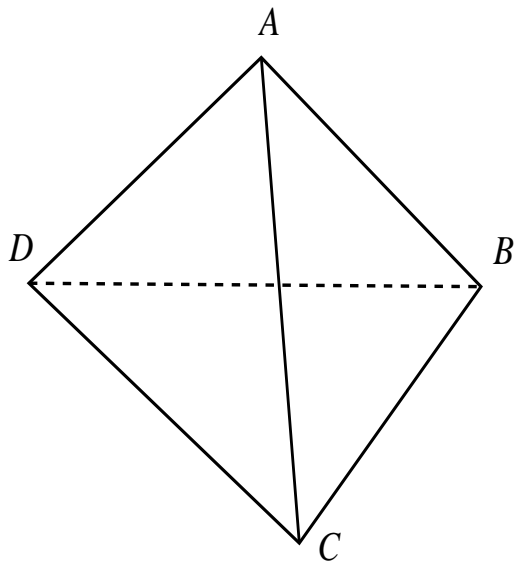}
\caption{A Whitehead move.}
\end{figure}

\thm{Definition}
{
The {triangulation graph ${\bf T}_N$} is a graph whose vertices are
isotopy classes of triangulations of $S^2$ on $N$ vertices, and there
is an edge joining 
nodes corresponding to $T_1$ and $T_2$ if and only if there exists a
Whitehead move transforming $T_1$ and $T_2$.
}

\thm{theorem}{\label{wconn}
The graph ${\bf T}_N$ is connected.
}

\begin{proof}
 Pick two distinguished vertices $v_1$ and $v_2$. 
For any starting triangulation $T$ there is a sequence of Whitehead
moves transforming $T$ into a triangulation $T_{v_1}$ where $v_1$ is
connected to every other vertex. Now consider the complement in
$T_{v_1}$ of $v_1$ and all the edges incident to it. This will be a
triangulation of an $(N-1)$-gon, all of whose vertices are those of the
$(N-1)$-gon. By a similar argument, this can be transformed by a
sequence of Whitehead moves into a triangulation where every vertex is
connected to $v_2$. Thus, it is seen that by a sequence of Whitehead
moves, every triangulation can be transformed to a particular
triangulation (where both $v_1$ and $v_2$ have valence $N-1$), and
thus ${\bf T}_N$ is connected. 
\end{proof}

\begin{proof}[of Theorem \ref{ctris}] It is enough to note that 
the cross ratio of any permutation $z_{\sigma(1)}, z_{\sigma(2)},
z_{\sigma(3)}, z_{\sigma(4)}$ is a rational function of the cross
ratio of $z_1, z_2, z_3, z_4$ (this is most easily seen when $z_1 = 0, z_2
= 1, z_3 = \infty$). By Lemma \ref{crossr} flipping the diagonal of
$ABCD$ corresponds to permuting the arguments of the cross ratio. Since
a permutation of the arguments corresponds to a fractional linear
transformation of the cross-ratio itself, Theorem \ref{ctris} follows.
\end{proof}

We summarize the results of this section for convenience:

\thm{theorem}{\label{teich}
$\mathcal{T}^N$ is a contractible $2 N - 6$ dimensional manifold, as
evidenced by coordinate systems coming from ideal triangulations of
$S_N^2$. Any two such coordinate systems are analytically equivalent. 
}

\section{Geometry of ideal tessellations}
\label{sec-cross}
First, let us review briefly the geometry of the upper half-space model
of $\hthree$. We will think of the ideal boundary $S_\infty^2$ of $\hthree$ as
the Riemann sphere $\cbar$. Hyperbolic planes are represented
by hemispheres whose equatorial circles are in $S_\infty^2$. In the
present context we think of straight lines as circles passing through
$\infty$. The corresponding hemispheres are vertical planes rising
above the lines.

Let $P$ be a convex polyhedron with all vertices on the ideal boundary
of $\hthree$. $P$ is the intersection of the half-spaces defined by its
faces. By an isometry of $\hthree$ and relabelling we can transform $P$ so
that the face $f_1$ of $P$ lies in the plane rising above the real
axis in $S_\infty^2$, and the vertices $v_1$, $v_2$ and $v_3$ are $0$,
$1$ and $\infty$ respectively. Furthermore, without loss of generality
we assume that $P$ lies above the half-plane $\Im(z)\geq0$.

The rest of the faces of $P$ are then oriented in such a way that the
interior of the corresponding hemispheres lie {\em outside} of $P$. 

$P$ defines a Euclidean tessellation of $\cbar$ in the natural
way: $P$ casts a shadow on the ideal boundary of $\hthree$ under the
orthogonal projection. The edges of $P$ are then mapped to
straight-line segments, and the faces of $P$ to convex polygons.
Denote the resulting tessellation of $\comps$ by $T_P$. This
tessellation has the following properties:

\medskip\noindent
{\bf Condition 1.} Every face $F$ of $T_P$ is inscribed in the circle
$\ccirc{F}$. 

\medskip\noindent
{\bf Condition 2.} No vertices of $T_P$ are contained in the
interior of $\ccirc{F}$. 

\medskip\noindent
{\bf Condition 3.} $T_P$ is contained in the upper half-plane of
$\comps$. 

\medskip\noindent
In the sequel we will assume for simplicity that $T_P$ is a
triangulation (unless otherwise indicated). Any more general
tessellation can be subdivided until it is a triangulation. First we
note the following: 

\thm{Lemma}
{\label{localT}
Condition 2 of is equivalent to the following:

\medskip\noindent
{\bf Condition $2'.$} for any two abutting triangles $ABC$ and $ADC$
of $T_P$, $D$ is not in the interior of $\ccirc{ABC}$.
}

\begin{proof}
Consider $P$. Lemma \ref{localT} is  equivalent to the
observation that the polyhedron $P$ is convex (Condition 2) if and
only if all of its edges are convexly bent (Condition $2'$).
\QED

\medskip\noindent
{\bf Note.} A simple direct Euclidean proof of Lemma \ref{localT} is
also possible. This will be left as an exercise for the reader.

\thm{Corollary}
{\label{WhiteMT}
Given an arbitrary triangulation $T'$ on the same vertex set as $T_P$,
$T'$ can be transformed into $T_P$ by a finite sequence of Whitehead
moves of the following kind: whenever $ABC$ and $ADC$ are abutting
triangles of $T'$ such that $D$ lies inside $\ccirc{ABC}$, we
change $ABC$ and $ADC$ into $ABD$ and $CBD$.
}

\begin{proof}
A Whitehead move of the described type corresponds to filling in a
missing tetrahedron $ABCD$ of a polyhedron lying above $T'$. Every
time a move as above happens, the edge $AC$ is buried, never to be
seen again. Since the number of possible edges is finite, the result
follows. 
\QED

The following facts from elementary Euclidean geometry will be needed
in the sequel. 

\thm{Lemma}{
\label{halfa}
Let $\mathcal{C}$ be a circle with center $O$ and let $ABC$ be a triangle
inscribed in $\mathcal{C}$. Then the  $\angle ACB=\half \angle AOB$ if
$C$ and $O$ are on the same side of $AB$ and $\angle ACB =\pi -\half
\angle AOB$ otherwise.
}

\begin{figure}
\label{douba}
\includegraphics*[width=3in,height=3in]{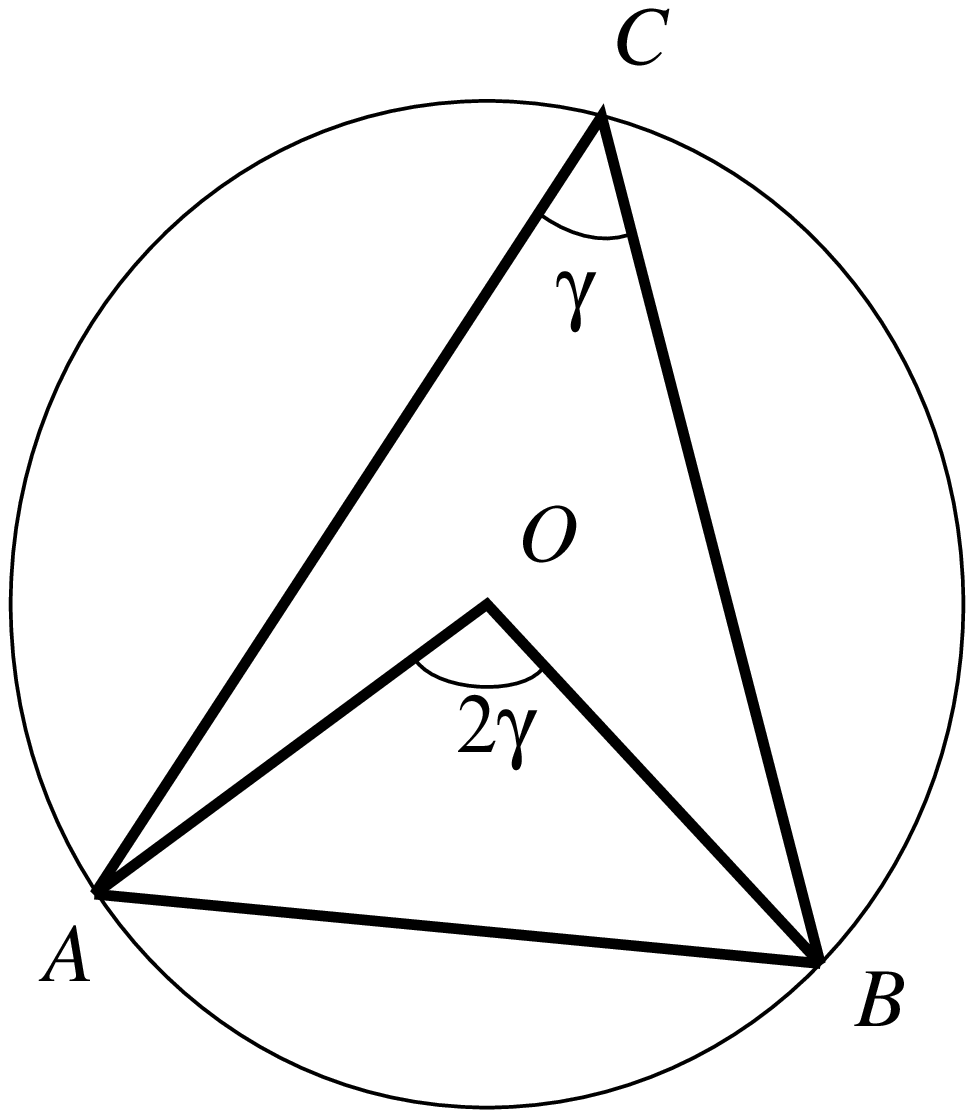}
\caption{}
\end{figure}

\thm{Lemma}{
\label{quad}
Let $ABCD$ be a quadrilateral. 

\begin{description}

\item{} $D$ is {\em outside} $\ccirc{ABC}$ if $\angle D + \angle B <\pi$.

\item{} $D$ is {\em inside} $\ccirc{ABC}$ if $\angle D + \angle B
>\pi$.

\item{} $D$ is {\em on} $\ccirc{ABC}$ if $\angle D + \angle B =\pi$.

\end{description}
}

The following important fact follows from Lemma \ref{halfa}:

\thm{theorem}{\label{diang}
The dihedral angle between the faces $ABC$ and $ADC$ of $P$ is equal
to the sum of angles $\angle ABC$ and $\angle ADC$.
}

\begin{proof}
First, observe that the angle between the the faces $ABC$
and $ADC$ is equal to the angle between the circles $\mathcal{C}_{ABC}$
and $\ccirc{ADC}$. The rest of the proof is contained in Figure 4.
\QED 

\begin{figure}
\label{dihed}
\includegraphics*[width=3in,height=3in]{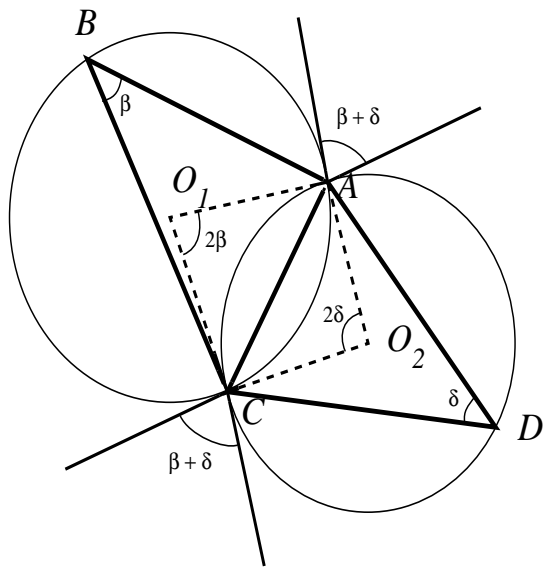}
\caption{Dihedral angle}
\end{figure}

The sum of the angles $\angle ABC$ and $\angle ADC$ is easily seen to
equal the argument of the cross-ratio

$$\crat(A, B, C, D)=[B, C, A, D]=\frac{(B-A)}{(B-C)} \frac{(C-D)}{(A-D)}$$

If $A$, $B$, $C$, and $D$ are transformed  by a hyperbolic isometry in
such a way that $A=\infty$, $B=1$, $C=0$, and $D=z$,
then $\crat(A, B, C, D) = z$, then Theorem \ref{shifts} below follows
from the discussion of section \ref{sec-intrinsic}.

\thm{theorem}{\label{shifts} 
With notation as above  $ \sh (ABC, ADC) = \log |c(A, B, C, D)|$. 
}

\begin{proof}
This is just a ``bent'' three-dimensional version of
Lemma \ref{crossr}. 
\end{proof}

The above observations allow us to prove Theorems \ref{mapcont} and
\ref{closed}, which are two of the steps of the proof of Theorem
\ref{charRI}.

\thm{theorem}{\label{mapcont}
The map $\mathfrak{g}$ is continuous.
}

\begin{proof}
The simpler case is one where $T_P$ is a genuine triangulation. Then,
it is clear that a small perturbation of the vertices of $T_P$ doesn't
change the combinatorics of $T_P$, and so continuity follows from
Theorem \ref{shifts} and the continuity of the cross-ratio. Things are
very slightly more complicated when $T_P$ has non-triangular faces.
Then, $T_P$ is combinatorially unstable: a small perturbation in the
vertices changes the combinatorial structure, but only in the
following simple fashion:

\thm{Lemma}{\label{charts}
For a sufficiently small $\epsilon$, a perturbation $T_{P^\epsilon}$ of
$T_P$ is combinatorially equivalent to $T_P$ with some diagonals added
to the non-triangular faces.
}

\begin{proof}[of Lemma] 
There exists an $\epsilon>0$, such that if a point $D$ is closer than
$\epsilon$ to $\ccirc{ABC}$ then $A$, $B$, $C$, and $D$ are
co-circular, for any triangle $ABC$ and vertex $D$ of $T_P$. 
\QED

Every way of adding diagonals to $T_P$ until we get a triangulation
corresponds to a different coordinate system of $T_N$ (where $N$ is
the number of vertices of $P$), and Lemma \ref{charts} shows that
every sufficiently small perturbation of $T_P$ is close to $T_P$ in at
least one of the coordinate systems. Since the transition maps between
the various coordinate systems are continuous (Theorem \ref{ctris}),
Theorem \ref{mapcont} follows. \QED

\thm{theorem}{\label{closed}
The image of $\mathfrak{g}$ is closed.
}

\begin{proof}
Let $\mathfrak{H}_1, \dots, \mathfrak{H}_k, \dots$ be a sequence
of metrics on $S_N^2$ converging to a metric $\mathfrak{H}$. Let $P_i$ be
such that 
$\mathfrak{g}(P_i) = \mathfrak{H}_i$. We will show that there exists a
$P_\infty$ such that 
$\mathfrak{g}(P_\infty) =\mathfrak{H}$. First choose a subsequence $P'_1, \dots,
P'_i, \dots$ such 
that all of the $P'_i$ have the same combinatorics. This is possible
since the number of possible combinatorial structures is finite. As
before, the vertices and faces of $P'_i$ are labelled in such a way
that $v_1(P'_i)=0$, $v_2(P'_i) = 1$, $v_3(P'_i)=\infty$, and
$f_1(P'_i)$ lies above the real axis. By compactness of the sphere
$\cbar$, there exists a limiting tessellation $T_{P'}$. If
$T_{P'}$ is non-degenerate (that is, no two vertices of a triangle have
coalesced into one), then by Theorem \ref{mapcont} it follows that
$\mathfrak{g}(P') = \mathfrak{H}$. 

We will show that $T_{P'}$ is always non-degenerate. If this is not
the case, let $t_i$ be a collapsing face of $T_{P'}$ which is abutting
a non-collapsing face $t_j$. Such a pair of faces must exist, since at
least one face ($f_1$) is not collapsing. By relabelling and a
hyperbolic isometry, send $t_j$ to the triangle $0, 1, \infty$,
and $t_i$ to $0, \infty, z$. Since $\sh(t_j, t_i) = \log |z|$, and
$\mathfrak{H}$ is a non-degenerate metric it follows that $z$ stays away
from $0$ and $\infty$, and so $t_i$ is not collapsing after all.
\QED

\thm{Remark}{
{\rm 
Theorem \ref{diang} and subsequent discussion is easily seen to
lead to the following pleasing ``hyperbolic'' interpretation of planar
triangulations. Consider a (not necessarily convex) polygon $Q$ in the
plane, such that the interior of $Q$ is triangulated in such a way
that edges of $Q$ are edges of the triangulation $T(Q)$. Then $T(Q)$
is the projection of an ideal polyhedron $\tilde{Q}$ onto the plane
$\comps$ at infinity of $H^3$, such that:

\begin{enumerate}
\item{} $\tilde{Q}$ has vertices at the vertices of $T(Q)$, plus one
vertex at the point $\infty$ of $\bar\comps.$

\item{} $Q$ is similar to the link of the vertex $v_\infty$ of
$\tilde{Q}$ at $\infty.$

\item{} $\tilde{Q}$ is star-shaped with respect to $v_\infty.$

\item{} If $v_1, \dots, v_k$ are the vertices of $Q$, then the
dihedral angle of $\tilde{Q}$ corresponding to the edge $v_k v_\infty$
is the Euclidean angle of $Q$ at $v_k.$

\item{} The dihedral angle of $\tilde{Q}$ corresponding to the
boundary edge $v_i v_{i+1}$ is equal to the euclidean angle at the
third vertex $w$ of the (unique) triangle $w v_i v_{i+1}$ of $T(Q)$
containing the edge $v_i v_{i+1}.$

\item{} The dihedral angle of $\tilde{Q}$ corresponding to a
non-boundary edge $AB$ of $T(Q)$ is equal to the sum of the angles at
$C$ and $D$ of the two triangles $ACB$ and $ADB$ abutting along the
edge $AB$

\end{enumerate}

The triangulation $T(Q)$ could also be taken to be {\em
immersed}, in which case all of the above statements still hold, with
the obvious changes in interpretation. 
}
}

\thm{Definition}{
A set of edges $C= \{e_1, \dots, e_k\}$ in a graph $G$ is called a {\em
cutset}, if the removal of those edges disconnects $G.$ A cutset $C$ is
called {\em minimal}, if no subset of $C$ is a cutset.
}

The simplest example of a cutset is the set of edges incident to a
single vertex of $G.$

The correspondence above can be used to prove the following result:

\thm{theorem}{\label{steiness}
Let $e_1, \dots, e_k$ be a minimal cutset of the 1-skeleton $\tilde
Q$ Then the sum $\Sigma$ of dihedral angles at $e_1, \dots, e_k$ is
strictly smaller than $(k-2) \pi$ if $e_1, \dots, e_k$ are not all
incident to one vertex. If $e_1, \dots, e_k$ are all incident to one
vertex then $\Sigma$ is exactly  $(k-2)\pi.$}

\begin{proof}
(also see figure 5) It is easy to see that
any minimal 
cutset as above is
actually the set of internal edges of a triangulation $T(\mathcal{A})$ of
an annulus $\mathcal{A}$ (possibly with one boundary component collapsed to a point $v$ ,
if all of the $e_i$ are incident to $v.$) From now on all references
will be to quantities in $T(Q)$ Let the inner and outer boundary
components of 
$\mathcal{A}$ be $\mathcal{A}_1$ and $\mathcal{A}_2$, respectively. The
edges of $T(\mathcal{A})$ naturally fall into three categories -- outer
boundary edges, inner boundary edges and internal edges. Similarly,
divide the angles of the triangles of $T(\mathcal{A})$ into the three
sets -- $A$ (angles opposite outer boundary), $B$ (angles opposite
inner boundary) and $\Gamma$ (angles opposite inner edges). Obviously, 
\begin{equation}
\label{abc}
\sum A + \sum B + \sum \Gamma = k \pi
\end{equation} 
(where $k$ is the number
of triangles and the cardinality of the cutset). Furthermore, by
Theorem \ref{diang}, 
\begin{equation}
\label{sg}
\sum \Gamma = \Sigma.
\end{equation}
 Now, note that the sum of
the angles incident (not opposite) to $\mathcal{A}_2$ is $(\card \mathcal{A}_2 - 2) \pi$, and further note that this sum is equal to $\sum B + 
(\card \mathcal{A})_2 \pi - \sum A.$ That is true since if $\alpha$ is {\em
opposite} to $\mathcal{A}_2$, then the other angles of the triangle
containing $\alpha$ are {\em incident} to $\mathcal{A}_2.$
Now, since 
\begin{equation}
\label{ab}
(\card \mathcal{A}_2 - 2) \pi = \sum B + \pi \card \mathcal{A}_2 - \sum A,
\end{equation}
it follows that $\sum A - \sum B = 2\pi.$ Since $\sum B$ is
greater than zero precisely when the inner boundary of $\mathcal{A}$ is
non-degenerate, it follows that $\sum A+\sum B > 2\pi$ whenever the
inner boundary of $\mathcal{A}$ is non-degenerate and $\sum A+\sum B =
2\pi$ otherwise. The statement of the theorem then follows from
equations \ref{abc}, \ref{sg}, and \ref{ab}.
\QED

\begin{figure}
\label{steiproof}
\includegraphics*[width=3in,height=3in]{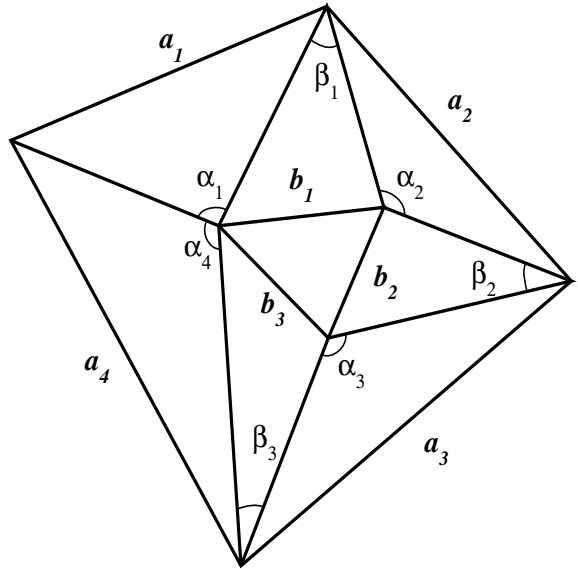}
\caption{Cutset sum}
\end{figure}

The above theorem is slightly stronger than Theorem 1 of
\cite{rivi90}, which is stated only for convex polyhedra. 
It turns out that the conditions of Theorem \ref{steiness} together
with the convexity conditions (dihedral angles are between $0$ and
$\pi$) completely characterize the sets of dihedral angles of convex
polyhedra. Proof of sufficiency is given in an upcoming paper
of the author.

\section{Ideal polyhedra are determined by their metric}

The purpose of this section is to prove Theorem \ref{Uniq}. First, a
couple of definitions:

\thm{Definition}{
\bd{(Generalized polyhedra and polygons)} \rm
A {\em generalized} convex hyperbolic 
polyhedron is represented in the projective model of $\hthree$ by a Euclidean
convex polyhedron which may have some vertices on or outside the
sphere at infinity (called ``infinite'' and
``hyperinfinite'' vertices respectively). However, each
edge must contain some points inside hyperbolic space.
 We will usually only be concerned with the part of a generalized
 polyhedron lying within $\hthree$.
Generalized convex polygons in $\mathbb{H}^2$ are defined similarly.
}

\thm{Definition}{\bd{(Links of vertices)} \rm
A generalized hyperbolic polyhedron has vertices of
three types: finite, hyperinfinite and infinite
vertices.

The ``link'' of a finite vertex of a
polyhedron is the spherical polygon obtained by
intersecting a small sphere centered at the vertex with
the polyhedron, and rescaling so the sphere has radius
$1$. So the edge lengths in the link are precisely the
face angles at the vertex.

For each hyperinfinite vertex there is a unique
hyperbolic plane orthogonal to the faces meeting at the
vertex.  The intersection of this plane with these faces
is a hyperbolic polygon which we will call the ``link'' of
the vertex. The edge lengths in the link are precisely the
lengths of common perpendiculars to adjacent sides
meeting at the hyperinfinite vertex.

For each infinite vertex there a 1-parameter
family of horospheres centered at the vertex. Each small
horosphere intersects the polyhedron in a Euclidean
polygon, which we will call the ``link'' of the vertex.
In this case the link is only well defined up to
Euclidean similarities.
}

\noindent{\bf Remark.} The link of an infinite or hyperinfinite vertex
then determines the corresponding {\em end} of the polyhedron up to
congruence.

\thm{Remark}{\label{linki}
{\rm
The link of an ideal vertex of $P$ is a Euclidean convex
polygon. Theorem \ref{shifts} shows that
if all vertices of $P$ are ideal, then the logarithm of the  ratio of
two adjacent sides of the link of a vertex $v$ is equal to the shear
between the two corresponding faces.
}}

The following result is obtained in \cite{rivi86} (see \cite{HR1},
Theorem 4.10):

\thm{theorem}{
\label{links}
A generalized convex polyhedron $P$ in hyperbolic $3$-space is
determined up to congruence by the type of its vertices and the edge
lengths of the links of its vertices.
}

\thm{Note}{
{\rm The edges of $P$ are not required to be non-degenerate, so some
of the dihedral angles may be $\pi$.
}}

This theorem means that two combinatorially equivalent polyhedra $P_1$
and $P_2$ such that the corresponding sides of corresponding links of
$P_1$ and $P_2$ are equal are congruent. 

\begin{proof}[of Theorem \ref{Uniq}] 
Let $M$ be a complete finite-volume hyperbolic surface homeomorphic to
$S_N^2$. Let $P_1$ and $P_2$ be two different embeddings of $M$ into
$\hthree$ as convex polyhedra. If $P_1(M)$ and $P_2(M)$ are
combinatorially equivalent, Theorem \ref{links} implies that $P_1(M)$ and
$P_2(M)$ are congruent. 

Assume that $P_1(M)$ and $P_2(M)$ are not combinatorially equivalent.
Then $P_1$ and $P_2$ induce two different cell decompositions $Q_1$
and $Q_2$ of $M$, where the edges of $Q_i$ are preimages of
corresponding edges of $P_i(M)$. Produce a new cell decomposition $Q$ of
$M$ by {\em superimposing} $Q_1$ and $Q_2$. The vertex set of $Q$ is
the union of $V$ (the cusps of $M$) with the set $V'$ of intersections
of edges of $Q_1$ with those of $Q_2$. The image of $Q$ under $P_i$
will be $P_i(M)$ with some extra edges and vertices drawn on it. Then
we can treat $P_1(M)$ and $P_2(M)$ as being of the same type (that of
$Q$), and Theorem \ref{links} may be applied. Thus $P_1(M)$ and
$P_2(M)$ are congruent. 
\QED

\section{Directions for further research}

To the author, the most painful shortcoming of the results presented
in this paper is the lack of any constructive method of producing
an embedding of a hyperbolic $N$-punctured sphere $\mathfrak{S}$ into
$\hthree$ as a convex ideal polyhedron. In particular, the tesselation
of $\mathfrak{S}$ induced by such an embedding is clearly canonical (in
view of Theorem \ref{Uniq}), and yet there seems no
known method of producing it.

The simple Theorem \ref{diang} turns out to be very useful. 
A number of consequences are given in the author's paper \cite{rvol}). 
An efficient
algorithm for producing a convex ideal polyhedron with prescribed dihedral
angles is contained in an upcoming joint paper of the author and
Warren~D.~Smith.

\paragraph{Acknowledgements.} The author would like to thank Craig
Hodgson and Warren~D.~Smith for valuable comments on earlier drafts of
this paper.


\end{document}